# The monoidal centre as a limit

Ross Street

**Abstract**

The centre of a monoidal category is a braided monoidal category. Monoidal categories are monoidal objects (or pseudomonoids) in the monoidal bicategory of categories. This paper provides a universal construction in a braided monoidal bicategory that produces a braided monoidal object from any monoidal object. Some properties and sufficient conditions for existence of the construction are examined.

## 1. Introduction

During question time after a talk [St2] at the Fields Institute, Peter Schauenburg asked whether the centre construction on a monoidal category (see [JS]) would fit into the general framework of [DMS] that I was describing. At the time I could not see how to do it. Reinforced by Peter's interest, the question stayed with me. During preparation of the paper [St3] on descent theory, intended for a publication arising from the same Fields Institute workshop, the answer began to dawn on me. Another topic at the top of my mind recently (in work with Michael Batanin and Alexei Davydov) has been Hochschild cohomology, and this too turns out to be relevant.

## 2. The centre of a monoidal object

In any monoidal bicategory $\mathcal{M}$, with tensor product $\otimes$ and unit $I$, we use the terms *pseudomonoid* and *monoidal object* for an object $A$ equipped with a binary multiplication $m : A \otimes A \longrightarrow A$ and a unit $j : I \longrightarrow A$ which are associative and unital up to coherent invertible 2-cells. A *monoidal morphism* $f : A \longrightarrow A'$ is a morphism equipped with coherent 2-cells $m \circ (f \otimes f) \Rightarrow f \circ m$ and $j \Rightarrow f \circ j$. The monoidal morphism is called *strong* when the coherent 2-cells are both invertible. A *monoidal 2-cell* is one compatible with these last coherent 2-cells. With the obvious compositions, this defines a bicategory Mon$\mathcal{M}$ of pseudomonoids in $\mathcal{M}$. For example, if $\mathcal{M}$ is the cartesian-monoidal 2-category Cat of categories, functors and natural transformations then Mon$\mathcal{M}$ is the 2-category MonCat of monoidal categories, monoidal functors and monoidal natural transformations as defined in [EK].

We now suppose $\mathcal{M}$ is braided. In fact, by the coherence result of [GPS], we suppose $\mathcal{M}$ is a braided Gray monoid in the sense of [DS]. The braiding for $\mathcal{M}$ is denoted by

$$c_{X,Y} : X \otimes Y \xrightarrow{\sim} Y \otimes X.$$

For any monoidal object $A$ of $\mathcal{M}$, a morphism $u : U \longrightarrow A$ is called a *centre piece* when it is equipped with an invertible 2-cell

$$\begin{array}{ccc} U \otimes A & \xrightarrow{c_{U,A}} & A \otimes U \\ {\scriptstyle u \otimes 1_A} \downarrow & \overset{\gamma}{\Rightarrow} & \downarrow {\scriptstyle 1_A \otimes u} \\ A \otimes A & & A \otimes A \\ & \searrow_m \quad A \quad \swarrow_m & \end{array}$$

such that the following equality holds.

[Diagram: large hexagonal pasting diagram involving $U \otimes A \otimes A \xrightarrow{c_{U,A\otimes A}} A \otimes A \otimes U$, with 2-cells $c_{1,m}$, $\gamma$, and associativity isomorphisms, equal to a second diagram involving $c_{U,A} \otimes 1$, $1 \otimes c_{U,A}$, $\gamma \otimes 1$, and $1 \otimes \gamma$.]

$\|$

A *morphism* $\sigma : u \Rightarrow v : U \longrightarrow A$ *of centre pieces* is a 2-cell $\sigma : u \Rightarrow v$ such that

$$\gamma \left( m \circ (\sigma \otimes 1_A) \right) = \left( m \circ (1_A \otimes \sigma) \circ c_{U,A} \right) \gamma .$$

We write $\mathcal{CP}(U, A)$ for the category of centre pieces so obtained. Using the pseudo-naturality of the braiding for $\mathcal{M}$, we see that we have a pseudofunctor

$$\mathcal{CP}(-, A) : \mathcal{M}^{\mathrm{op}} \longrightarrow \mathrm{Cat}$$

defined on morphisms $f : V \longrightarrow U$ by composition; that is, the functor

$$\mathcal{CP}(f, A): \mathcal{CP}(U, A) \longrightarrow \mathcal{CP}(V, A)$$

takes a centre piece $u$ with $\gamma$ to $u \circ f$ with the 2-cell obtained by pasting the square containing the pseudonaturality isomorphism $c_{f, 1_A}$ onto the top of the pentagon containing $\gamma$.

The *centre of* $A$ is a birepresenting object $\mathcal{Z}A$ (in the sense of [St0]) for the pseudofunctor $\mathcal{CP}(-, A)$. This means we have a centre piece $i: \mathcal{Z}A \longrightarrow A$, composition with which induces an equivalence of categories

$$\mathcal{M}(U, \mathcal{Z}A) \simeq \mathcal{CP}(U, A).$$

It follows that the centre of $A$ is unique up to equivalence if it exists.

**Proposition 2.1** *The centre $\mathcal{Z}A$ of a monoidal object $A$ is a braided monoidal object in the sense of* [DS], *and the morphism* $i: \mathcal{Z}A \longrightarrow A$ *is strong monoidal.*

**Proof** The composite $\mathcal{Z}A \otimes \mathcal{Z}A \xrightarrow{i \otimes i} A \otimes A \xrightarrow{m} A$ equipped with the 2-cell

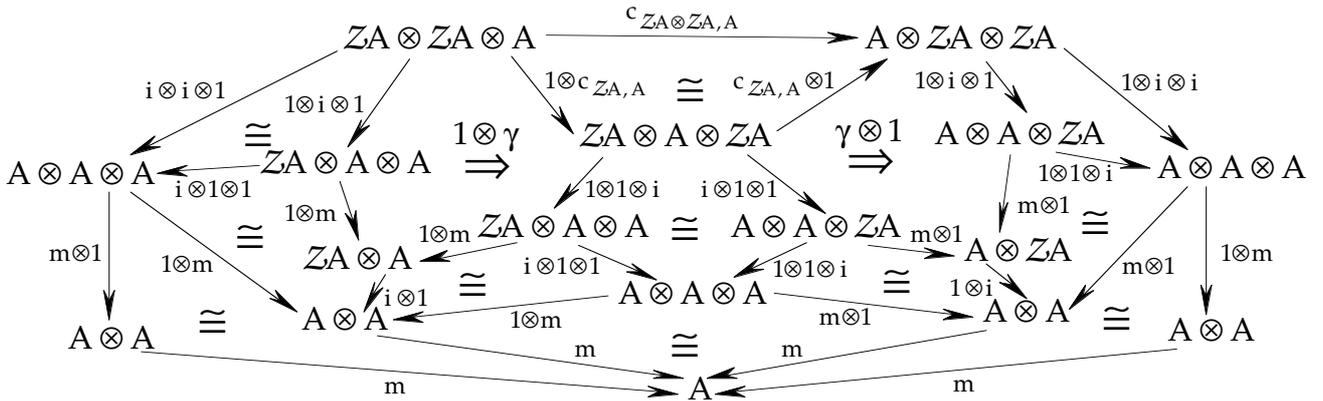

is a centre piece. So, up to a unique invertible 2-cell, there is a morphism

$$m: \mathcal{Z}A \otimes \mathcal{Z}A \longrightarrow \mathcal{Z}A$$

and an invertible 2-cell $i \circ m \cong m \circ (i \otimes i)$ compatible with the 2-cells of the two centre pieces. Also $j: I \longrightarrow A$, equipped with the obvious 2-cell, is a centre piece and so induces a morphism $j: I \longrightarrow \mathcal{Z}A$ with $i \circ j \cong j$. Largish pasting diagrams prove that $\mathcal{Z}A$ becomes monoidal with $i: \mathcal{Z}A \longrightarrow A$ strong monoidal. The braiding for $\mathcal{Z}A$ is the invertible 2-cell

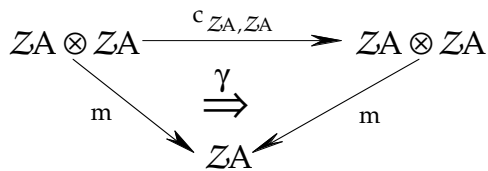

whose composite with $i: \mathcal{Z}A \longrightarrow A$ is the pasting composite below.

$$\begin{array}{c}
\mathcal{Z}A \otimes \mathcal{Z}A \xrightarrow{c_{\mathcal{Z}A, \mathcal{Z}A}} \mathcal{Z}A \otimes \mathcal{Z}A \\
{}_{1\otimes i}\swarrow \quad \overset{c_{1,i}}{\cong} \quad \searrow{}_{i\otimes 1} \\
\mathcal{Z}A \otimes A \xrightarrow{c_{\mathcal{Z}A, A}} A \otimes \mathcal{Z}A \\
{}_{i\otimes 1}\downarrow \quad \overset{\gamma}{\Rightarrow} \quad \downarrow{}_{1\otimes i} \\
A \otimes A \qquad A \otimes A
\end{array}$$

with outer arrows $m$ down to $\mathcal{Z}A$ on each side, $\cong$, then $i$ down to $A$, and inner $m$ arrows to $A$.

The braiding axioms follow from the defining property of a centre piece. **Q.E.D.**

The following two propositions are routinely proved.

**Proposition 2.2** *If* $F : \mathcal{M} \longrightarrow \mathcal{N}$ *is a braided monoidal pseudofunctor and* $u : U \longrightarrow A$ *is a centre piece for the monoidal object* $A$ *in* $\mathcal{M}$ *then* $Fu : FU \longrightarrow FA$ *is canonically a centre piece for the monoidal object* $FA$ *in* $\mathcal{N}$. *There is a canonical comparison morphism* $F\mathcal{Z}A \longrightarrow \mathcal{Z}FA$ *provided the centres of* $A$ *and* $FA$ *exist.*

**Proposition 2.3** *The pseudofunctor* $\mathcal{CP}(-, A) : \mathcal{M}^{op} \longrightarrow \mathrm{Cat}$ *preserves all weighted bicategorical limits that exist in* $\mathcal{M}^{op}$ *and, that as colimits in* $\mathcal{M}$, *are preserved by* $- \otimes A$.

## 3. Existence

In this section we shall look at conditions on the braided monoidal bicategory $\mathcal{M}$ for monoidal centres to exist. Because of Proposition 2.3, we expect $\mathcal{CP}(-, A) : \mathcal{M}^{op} \longrightarrow \mathrm{Cat}$ to be birepresentable when each $- \otimes A$ preserves colimits and a "special birepresentability theorem" applies to $\mathcal{M}$.

Recall from [DS] that a monoidal bicategory $\mathcal{M}$ is called *left [right] closed* when, for each object $B$, the pseudofunctor $- \otimes B : \mathcal{M} \longrightarrow \mathcal{M}$ [$B \otimes - : \mathcal{M} \longrightarrow \mathcal{M}$] has a right biadjoint (and so preserves bicategorical colimits). We call $\mathcal{M}$ closed when it is both left and right closed; we denote the right biadjoint of $- \otimes B$ by $[B, -] : \mathcal{M} \longrightarrow \mathcal{M}$ and we have a family of equivalences

$$\mathcal{M}(A \otimes B, C) \simeq \mathcal{M}(A, [B, C]),$$

pseudonatural in each variable, and called the *closedness equivalences*. Taking $A = [B, C]$,

we find an *evaluation morphism* ev : $[B,C] \otimes B \longrightarrow C$ that, up to isomorphism, is taken to the identity by the closedness equivalence.

Assume $\mathcal{M}$ is braided and left closed. It follows that $\mathcal{M}$ is closed. From any monoidal object A we shall construct a Hochschild-like truncated pseudo-cosimplicial object $C$A:

$$A \underset{\partial_1}{\overset{\partial_0}{\rightrightarrows}} [A,A] \underset{\partial_2}{\overset{\partial_0}{\underset{\partial_1}{\rightrightarrows}}} [A \otimes A, A]$$

as follows. The morphism $\partial_0 : A \longrightarrow [A,A]$ corresponds under the closedness equivalence to $m : A \otimes A \longrightarrow A$. The morphism $\partial_1 : A \longrightarrow [A,A]$ corresponds under the closedness equivalence to the composite $A \otimes A \xrightarrow{c_{A,A}} A \otimes A \xrightarrow{m} A$. The morphisms

$$\partial_0, \partial_1, \partial_2 : [A,A] \longrightarrow [A \otimes A, A]$$

correspond under the closedness equivalence to the morphisms

$$[A,A] \otimes A \otimes A \xrightarrow[\underset{(1_A \otimes \mathrm{ev}) \circ (c_{[A,A],A} \otimes 1_A)}{\mathrm{ev} \circ (1 \otimes m)}]{m \circ (\mathrm{ev} \otimes 1_A)} A .$$

One easily finds the coherent invertible 2-cells

$$\partial_0 \circ \partial_0 \cong \partial_1 \circ \partial_0 , \quad \partial_0 \circ \partial_1 \cong \partial_2 \circ \partial_0 , \quad \partial_2 \circ \partial_1 \cong \partial_1 \circ \partial_1 .$$

**Proposition 3.1** *The bicategorical limit of the pseudo-cosimplicial diagram $C$A is the centre of* A.

**Proof** The proof is by transport across the closedness equivalences. **Q.E.D.**

It is shown in [St0] how to construct this pseudo-descent-like limit in a bicategory $\mathcal{M}$ that admits finite products, iso-inserters, and cotensoring with the arrow category **2**.

**Corollary 3.2** *In any finitely complete, closed, braided monoidal bicategory, every monoidal object has a centre. Any braided monoidal pseudofunctor that is strong closed and finite-limit preserving preserves centres.*

Examples of such bicategories abound. Let $\mathcal{A}$ be any (small) braided promonoidal 2-category (in particular, $\mathcal{A}$ could be a braided monoidal category). Take $\mathcal{M}$ to be the 2-category $[\mathcal{A}, \mathrm{Cat}]$ of 2-functors from $\mathcal{A}$ to Cat. This is a complete and cocomplete 2-category and so is also as a bicategory. It becomes closed monoidal under the Day convolution tensor product defined by the coends

$$(F \otimes G)A = \int^{B,C} P(B,C;A) \times FB \times GC$$

in Cat. The braiding is induced by that on $\mathcal{A}$ and the symmetry on Cat. A good example of an appropriate $\mathcal{A}$ is provided by the category of automorphisms of a groupoid $\mathcal{G}$ as described in Section 7, Example 9 of [DS].

Alternatively we could take $\mathcal{M}$ to be the 2-category $\text{Hom}(\mathcal{A}, \text{Cat})$ of pseudofunctors from $\mathcal{A}$ to Cat. This is complete and cocomplete as a bicategory. It becomes a closed monoidal bicategory under the convolution tensor product defined by the pseudocoends

$$(F \otimes G)A = \int^{B,C}_{ps} P(B,C;A) \times FB \times GC$$

in Cat. Again, the braiding is induced by that on $\mathcal{A}$ and the symmetry on Cat.

# References

[DMS]  B. Day, P. McCrudden and R. Street, Dualizations and antipodes, *Applied Categorical Structures* (to appear).

[DS]  B. Day and R. Street, Monoidal bicategories and Hopf algebroids, *Advances in Math.* **129** (1997) 99-157.

[EK]  S. Eilenberg and G.M. Kelly, Closed categories, *Proc. Conf. Categorical Algebra at La Jolla 1965* (Springer-Verlag, Berlin 1966) 421-562.

[GPS]  R. Gordon, A.J. Power and R. Street, Coherence for tricategories, *Memoirs Amer. Math. Soc.* **117** (1995) #558.

[JS]  A. Joyal and R. Street, Tortile Yang-Baxter operators in tensor categories, *J. Pure Appl. Algebra* **71** (1991) 43-51.

[St0]  R. Street, Fibrations in bicategories, *Cahiers topologie et géométrie différentielle* **21** (1980) 111-160; **28** (1987) 53-56.

[St1]  R. Street, The quantum double and related constructions, *J. Pure Appl. Algebra* **132** (1998) 195-206.

[St2]  R. Street, Formal representation theory, Talk at the "Workshop on Categorical Structures for Descent and Galois Theory, Hopf Algebras and Semiabelian Categories" (Fields Institute, 24 September, 2002); slides and audio at <http://www.fields.utoronto.ca/audio/02-03/galois_and_hopf/street/>.

[St3]  R. Street, Categorical and combinatorial aspects of descent theory (submitted, March 2003) math.CT/0303175.